\documentclass[aps,pra, twocolumn, groupedaddress, nofootinbib]{revtex4}
\usepackage{amsmath,amsfonts,amssymb,amsthm, bbm}

\usepackage{graphicx}   % need for figures
\usepackage{subfigure}  % and subfigures

\newcommand\doi[1]{\href{http://dx.doi.org/#1}{doi:#1}}
\newcommand\arxivref[1]{\href{http://arxiv.org/abs/#1}{arXiv:#1}}

\newcommand{\pML}{\hat{p}_{\mathrm{ML}}}
\newcommand{\pLI}{\hat{p}_{\mathrm{LI}}}
\newcommand{\cL}{\mathcal{L}}
\newcommand{\diff}{\mathrm{d}\!}
\newcommand{\Ravg}{\overline{R}}

\usepackage{hyperref} % must be loaded last
\graphicspath{{./}{./figures/}}
\begin{document}

%---------------------------------------------------------------------------------------------------------------------
\title{Estimating the bias of a noisy coin}
\author{Christopher Ferrie}
\affiliation{
Institute for Quantum Computing,
University of Waterloo,
Waterloo, Ontario, Canada, N2L 3G1}
\affiliation{
Department of Applied Mathematics,
University of Waterloo,
Waterloo, Ontario, Canada, N2L 3G1}
\author{Robin Blume-Kohout}
\affiliation{
T-4 and CNLS, Los Alamos National Laboratory, Los Alamos NM 87545}

\date{\today}
\begin{abstract}
Optimal estimation of a coin's bias using \emph{noisy} data is surprisingly different from the same problem with noiseless data.  We study this problem using entropy risk to quantify estimators' accuracy.  We generalize the ``add $\beta$'' estimators that work well for noiseless coins, and we find that these \emph{hedged maximum-likelihood} (HML) estimators achieve a worst-case risk of $O(N^{-1/2})$ on noisy coins, in contrast to $O(N^{-1})$ in the noiseless case.  We demonstrate that this increased risk is unavoidable and intrinsic to noisy coins, by constructing minimax estimators (numerically).  However, minimax estimators introduce extreme bias in return for slight improvements in the worst-case risk.  So we introduce a \emph{pointwise} lower bound on the minimum achievable risk as an alternative to the minimax criterion, and use this bound to show that HML estimators are pretty good.  We conclude with a survey of scientific applications of the noisy coin model in social science, physical science, and quantum information science.
\end{abstract}
%\classification{}

\maketitle

%---------------------------------------------------------------------------------------------------------------------

%\tableofcontents

%---------------------------------------------------------------------------------------------------------------------
\section{Introduction}
%---------------------------------------------------------------------------------------------------------------------

Coins and dice lie at the foundation of probability theory, and estimating the probabilities associated with a biased coin or die is one of the oldest and most-studied problems in statistical inference.  But what if the data are \emph{noisy} -- e.g., each outcome is randomized with some small (but known) probability before it can be recorded?  Estimating the underlying probabilities associated with the ``noisy coin'' or ``noisy die'' has attracted little attention to date.  This is odd and unfortunate, given the range of scenarios to which it applies.  Noisy data appear in social science (in the context of randomized response), in particle physics (because of background counts), in quantum information science (because quantum states can only be sampled indirectly), and in other scientific contexts too numerous to list.

Noisy count data of this sort are usually dealt with in an ad hoc fashion.  For instance, if we are estimating the probability $p$ that a biased coin comes up ``heads'', and each observation gets flipped with probability $\alpha$, then after $N$ flips, the expectation of the number of ``heads'' observed ($n$) is not $Np$ but $Nq$, where
$$q = \alpha + p(1-2\alpha).$$
The ad hoc solution is to simply estimate $\hat{q} = \frac{n}{N}$ and invert $q(p)$ to get a \emph{linear inversion} estimator
\begin{equation}
\pLI = \frac{\hat{q}-\alpha}{1-2\alpha}. \label{eq:pLI}
\end{equation}
Rather awkwardly, though, $\pLI$ may be negative!  Constraining it to the interval $[0,1]$ fixes this problem, but now $\hat{p}$ is a biased estimator.  Whether it is a \emph{good} estimator -- i.e., accurate and risk-minimizing in some sense -- becomes hopelessly unclear.

In this paper, we analyze point estimators for noisy binomial data from the ground up, using \emph{entropy risk} as our benchmark of accuracy.  Most of our results for coins extend to noisy dice (multinomial data) as well.  We begin with a review of optimal estimators for noiseless coins, then examine the differences between noisy and noiseless data from a statistically rigorous perspective.  We show how to generalize good \emph{noiseless} estimators to noisy data, compare their performance to (numerical) minimax estimators for noisy coins, point out the shortcomings in the minimax approach, and propose an alternative estimator with good performance across the board.

\section{Starting points: Likelihood and Risk}

The noisy coin is a classic point estimation problem:  we observe data $n$, sampled from one of a parametric family of distributions $Pr(n|p)$ parameterized by $p$, and we seek a point estimator $\hat{p}(n)$ that minimizes a risk function $R(\hat{p};p)$.  The problem is completely defined by (1) the sampling distribution, and (2) the risk function.

\subsection{The sampling distribution and likelihood function}

For a \emph{noiseless coin}, the sampling distribution is
\begin{equation*}
Pr(n|p) = \binom{N}{n}p^n(1-p)^{N-n},
\end{equation*}
where $N$ is the total number of samples, and $p\in[0,1]$.  Our interest, however, is in the \emph{noisy coin}, wherein each of the $N$ observations gets flipped with probability $\alpha \in \left[0,\frac12\right)$.  The sampling distribution is therefore
\begin{equation*}
Pr(n|p) = \binom{N}{n}q^n(1-q)^{N-n},
\end{equation*}
where the ``effective probability'' of observing heads is
\begin{equation*}
q(p) = \alpha + p(1-2\alpha) = p + \alpha\left(1-2p\right).
\end{equation*}
The noiseless coin is simply the special case where $\alpha=0$.

Everything that the data $n$ imply about the underlying parameter $p$ (or $q$) is conveyed by the likelihood function,
\begin{equation*}
\cL(p) = Pr(n|p).
\end{equation*}
$\cL$ measures the \emph{relative} plausibility of different parameter values.  Since its absolute value is never relevant, we ignore constant factors and use
\begin{equation*}
\cL(p) = \left(\alpha + p(1-2\alpha)\right)^n\left(1-\alpha - p(1-2\alpha)\right)^{N-n}
\end{equation*}
Frequentist and Bayesian analysis differ in how $\cL(p)$ should be used.  Frequentist estimators, including maximum likelihood (ML), generally make use \emph{only} of $\cL(p)$.  Bayesian estimators use Bayes' Rule to combine $\cL(p)$ with a prior distribution $P_0(p)\diff p$ and obtain a posterior distribution,
\begin{equation*}
Pr(p|n) = \frac{\cL(p)P_0(p)\diff p}{\int{\cL(p)P_0(p)\diff p}},
\end{equation*}
which determines the estimate.

\subsection{Risk and cost}

Since there is a linear relationship between $q$ and $p$, we could easily reformulate the problem to seek an estimator $\hat{q}$ of $q\in\left[\alpha,1-\alpha\right]$.  Estimation of $q$ and of $p$ are equivalent -- \emph{except} for the cost function.  If  for practical reasons we care about $q$, then the cost function will naturally depend on $q$.  But if we are fundamentally interested in $p$, then the cost function will naturally depend on $p$.  This determines whether we should analyze the problem in terms of $\hat{q}$ or $\hat{p}$.  It also implies (somewhat counterintuitively) that optimal estimators of $p$ and $q$ need not share the linear relationship between the variables themselves.

\emph{Entropy risk} provides a good illustration.  Suppose that the estimator $\hat{p}(n)$ [or $\hat{q}(n)$] is used to predict the next event.  The cost of event $k$ (e.g. ``heads'') occurring depends on the estimated $\hat{p}_k$, and is given by  $-\log\hat{p}_k$.  (This ``log loss'' rule models a variety of concrete problems, most notably gambling and investment).  The \emph{expected cost} depends on the true probability and on the estimate, and equals $C = -\sum_k{p_k\log\hat{p}_k}$.  Some of this cost is unavoidable, for random events cannot be predicted perfectly.  The \emph{minimum} cost, achieved when $\hat{p}_k = p_k$, is $C_{\mathrm{min}} = -\sum_k{p_k\log p_k}$.  The extra cost, which stems entirely from inaccuracy in $\hat{p}_k$, is the [entropy] \emph{risk}:
\begin{equation*}
R(\hat{p};p) = C - C_{\mathrm{min}} = \sum_k{p_k\left(\log p_k - \log \hat{p}_k\right)}.
\end{equation*}
This celebrated quantity, the \emph{relative entropy} or \emph{Kullback-Leibler divergence} \cite{Cover2006Elements} from $p$ to $\hat{p}$, is a good measure of how inaccurately $\hat{p}$ predicts $p$.

Now (returning to the noisy coin), both $q$ and $p$ are probabilities, so we could minimize the entropy risk from $q$ to $\hat{q}$, or from $p$ to $\hat{p}$.  But in the problems we consider, the events of genuine interest are the underlying (hidden) ones -- \emph{not} the ones we get to observe.  So $p$ is the operationally relevant probability, not $q$.  This is quite important, for if $p\approx0$ then $q\approx\alpha$.  Suppose the estimate is off by $\epsilon \ll \alpha$.  The \emph{risk} of that error depends critically on whether we compare the $q$'s,
$$R(\hat{q};q) \approx R\left(\alpha+\epsilon; \alpha\right) \approx \frac{\epsilon^2}{2\alpha(1-\alpha)},$$
or the $p$'s, in which case $R(\hat{p};p) = \infty$ if $\hat{p}=0$ and $p\neq0$, or (otherwise),
$$R(\hat{p};p) \approx R(\epsilon; 0) \approx \epsilon.$$
Entropy risk has special behavior near the state-set boundary -- i.e., when $p\approx0$ -- which has a powerful effect on estimator performance.

\section{The noiseless coin}

The standard biased coin, with bias $p = $ Pr(``heads''), is a mainstay of probability and statistics.  It provides excellent simple examples of ML and Bayesian estimation.

\subsection{The basics}

When a coin is flipped $N$ times and $n$ ``heads'' are observed, the likelihood function is
$$\cL(p) = p^n(1-p)^{N-n},$$
and its maximum is achieved at
$$\pML = \frac{n}{N}.$$
So maximum likelihood estimation agrees with the obvious na\"ive method of \emph{linear inversion} (Eq. \ref{eq:pLI}), which equates probabilities with observed frequencies.  Risk plays no role in deriving the ML estimator, and one objection to ML is that since $Pr(n=0|p)$ is nonzero for $p>0$, it is quite possible to assign $\pML=0$ when $p\neq0$, which results in \emph{infinite} entropy risk.  The practical problem here is that $\hat{p}=0$ may be interpreted as a willingness to bet at infinite odds against ``heads'' coming up -- which is more or less obviously a bad idea.

The simplest Bayesian estimator is \emph{Laplace's Law}\cite{Ristad1995Natural}.  It results from choosing a Lebesgue (``flat'') prior $P_0(p)\diff p = \diff p$, and then reporting the mean of the posterior distribution:
\begin{eqnarray}
\hat{p} &=& \int{pPr(p|n)\diff p} \nonumber \\
&=& \frac{\int{p\cL(p)P_0(p)\diff p}}{\int{\cL(p)P_0(p)\diff p}} \nonumber \\
&=& \frac{\int{p^{n+1}(1-p)^{N-n}\diff p}}{\int{p^{n}(1-p)^{N-n}\diff p}} \nonumber \\
&=& \frac{n+1}{N+2}.\nonumber
\end{eqnarray}
This estimator is also known as the ``add 1'' rule, for it is equivalent to (i) adding 1 fictitious observation of each possibility (heads and tails, in this case), then (ii) applying linear inversion or ML.

The derivation of Laplace's Law poses two obvious questions.  Why report the posterior mean? and why use the Lebesgue prior $\diff p$?  The first has a good answer, while the second does not.

We report the posterior mean because it minimizes the expected entropy risk -- i.e., it is the \emph{Bayes estimator}\cite{Lehmann1998Theory} for $P_0(p)\diff p$.  To prove this, note that in a Bayesian framework, we have assumed that $p$ is in fact a random variable distributed according to $P_0(p)\diff p$, and therefore \emph{given} the data $n$, $p$ is distributed according to $Pr(p|n)$.  A simple calculation shows that the expected entropy risk,
$$\Ravg = \int{R(\hat{p};p)Pr(p)\diff p},$$
is minimized by setting $\hat{p} = \int{p Pr(p) \diff p}$.  The posterior mean is Bayes for a large and important class of risk functions called \emph{Bregman divergences}, the most prominent of which is Kullback-Leibler divergence.  We will make extensive use of this convenient property throughout this paper, but in a broader context it is important to remember that for many other risk functions, the Bayes estimator is \emph{not} the posterior mean.

For the Lebesgue prior, on the other hand, there is little justification.  Lebesgue measure is defined on the reals by invoking translational symmetry, which does not exist for the probability simplex.  In fact, convenience is the best argument -- it certainly makes calculation easy!  But this argument can be extended to a large class of \emph{conjugate priors}.  Updating a Lebesgue prior in light of binomial (Bernoulli) data always yields a \emph{Beta distribution} \cite{betaref} as the posterior,
$$Pr(p) \propto p^{\beta-1}(1-p)^{\gamma-1}\diff p.$$
The family of Beta distributions is closed under the operation ``add another observation'', which multiplies $\cL(p)$ by $p$ or $(1-p)$.  So the same posterior form is obtained whenever the prior is a Beta distribution.  The conjugate priors for binomial data are therefore Beta distributions.

Any Beta prior is a convenient choice, but only those that maintain symmetry between heads and tails can be considered ``noninformative'' priors, e.g.
\begin{equation*}
P_0(p) \propto p^{\beta-1}(1-p)^{\beta-1},
\end{equation*}
for any real $\beta>0$.  The Bayes estimator for such a prior is (via Bayes Rule and the posterior mean),
\begin{equation*}
\hat{p} = \frac{n + \beta}{N+2\beta},
\end{equation*}
a rule known (for obvious reasons) as ``add $\beta$'', or (more obscurely) as \emph{Lidstone's Law}\cite{Ristad1995Natural}.

The ``add $\beta$'' estimators (including Laplace's Law) are convenient and sensible generalizations of the linear inversion estimator.  As we shall see later, they also generalize the ML estimator in an elegant way.  Better yet, they are Bayes estimators for Beta priors (with respect to entropy risk).  However, there are infinitely many \emph{other} priors, with their own unique Bayes estimators.

\subsection{Minimax and the noiseless coin}

For a true Bayesian, the previous section's analysis stands alone.  Every scientist or statistician must look into his or her heart, find the subjective prior therein, and implement its Bayes estimator.  For less committed Bayesians, the \emph{minimax} criterion\cite{Lehmann1998Theory} can be used to select a unique optimal estimator.  The essence of minimax reasoning is to compare different estimators not by their \emph{average} performance (over $p$), but by their \emph{worst-case} performance.  This eliminates the need to choose a measure (prior) over which to average.  Since this is a fundamentally frequentist notion, it is all the more remarkable that the minimax estimator is -- for a broad range of problems -- actually the Bayes estimator for a particular prior!

This theorem is known as \emph{minimax-Bayes duality}\cite{Lehmann1998Theory}.  The minimax estimator is, by definition, the estimator $\hat{p}(\cdot)$ that minimizes the maximum risk (thus the name):
\begin{equation*}
R_{\mathrm{max}}\left[\hat{p}(\cdot)\right] = \max_p{ \sum_n{Pr(n|p)R(\hat{p}(n);p)} }
\end{equation*}
Minimax-Bayes duality states that, as long as certain convexity conditions are satisfied (as they are for entropy risk), the minimax estimator is the Bayes estimator for a \emph{least favorable prior} (LFP), $P_\mathrm{worst}(p)\diff p$.

This duality has some useful corollaries.  The Bayes risk of the LFP (the minimum achievable average risk, which is achieved by the LFP's Bayes estimator) is equal to the minimax risk.  So the pointwise risk $R(\hat{p}_{\mathrm{Bayes}}(n);p)$ is identically equal to $R_{\mathrm{max}}$ at every support point of the LFP.  Furthermore, the LFP has the highest Bayes risk of all priors.  This means that the Bayes estimator of \emph{any} prior $P_0(p)\diff p$ yields both upper and lower bounds on the minimax risk.
\begin{equation} \label{eq:RiskBounds}
R_{\mathrm{Bayes}}(P_0) \leq R_{\mathrm{minimax}} \leq R_{\mathrm{max}}(P_0)
\end{equation}
These bounds coincide only for least favorable priors and minimax estimators.

Quite a lot is known about minimax estimators for the noiseless coin.  For small $N$, minimax estimators have been found numerically using the ``Box simplex method'' \cite{Wieczorkowski1998Calculating,Wieczorkowski1992Minimax}.  For large $N$, the binomial distribution of $n$ is well approximated by a Poisson distribution with parameter $\lambda=Np$.  By applying this approximation, and lower-bounding the maximum risk by the Bayes risk of the ``add 1'' estimator, it has been shown \cite{Rukhin1993Minimax}
\[
R_{\mathrm{max}}\left[\hat{p}(\cdot)\right] \geq   \frac1{2N}+ O\left(\frac1{N^2}\right).
\]
This limit can be achieved via the peculiar but simple estimator \cite{Braess2004Bernstein,Braess2002How}
\[
\hat p(n):=\begin{cases}
\frac {n+1/2}{N+5/4} & \textrm{if } n=0,\\
\frac {n+1}{N+7/4}   & \textrm{if } n=1,\\
\frac {n+3/4}{N+7/4} & \textrm{if } n=N-1,\\
\frac {n+3/4}{N+5/4} & \textrm{if } n=N,\\
\frac {n+3/4}{N+3/2} & \textrm{otherwise}.
\end{cases}
\]
This estimator is very nearly ``add $3/4$'', and in fact the ``add $\beta$'' estimators are almost minimax.  ``Add $1/2$'', which corresponds to the celebrated \emph{Jeffreys' Prior}\cite{Jeffreys1939Theory}, is a de facto standard.  However, it has been shown \cite{Rukhin1993Minimax,Krichevskiy1998Laplaces} that the very best ``add $\beta$'' estimator (although not \emph{quite} minimax) is the one with $\beta = \beta_0 = 0.509\ldots$, whose asymptotic risk is
\[
R_{\mathrm{max}}\left[\hat{p}(\cdot)\right] ={\beta_0}\frac1N+ O\left(\frac1{N^2}\right).
\]
These estimators \emph{hedge} effectively against unobserved events (in contrast to ML, which assigns $p=0$ if $n=0$), and we will strive to generalize them for the noisy coin.

\section{The noisy coin}

\begin{figure}[h]
%\resizebox{.7\columnwidth}{!}
\begin{tabular}{c}
\includegraphics[width=2in]{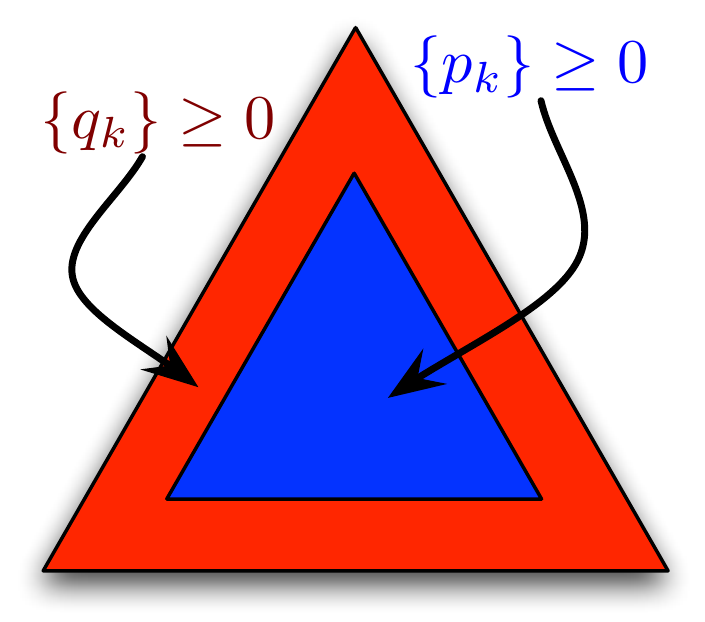} \\ \includegraphics[width=3in]{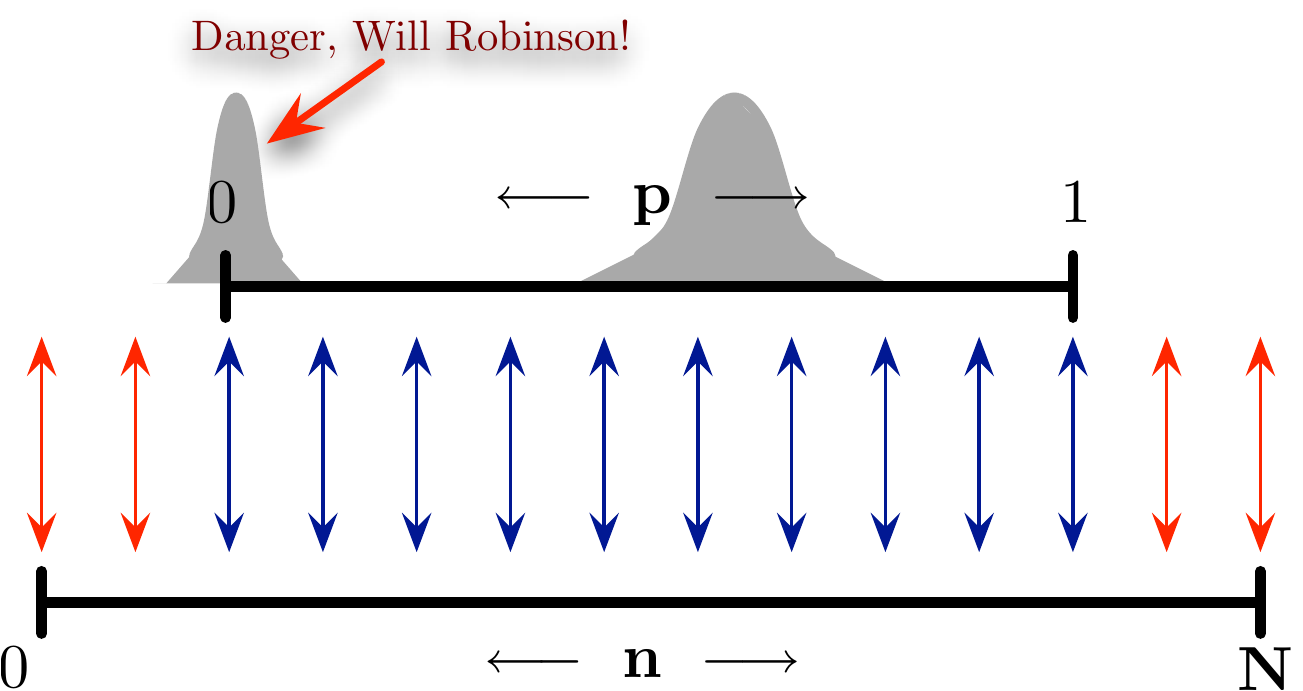}
\end{tabular}
\caption{Most of the complications in noisy coin estimation come from the existence of -- and linear relationship between --\emph{two} simplices, one containing all the true probability distributions $\{p_k\}$, and the other containing all effective distributions $\{q_k\}$.  Above (top), this is illustrated for a 3-sided die (not considered in this paper).  Below, an annotated diagram of \emph{sampling mismatch} for a coin shows how standard errors for the linear inversion estimator can easily extend outside the $p$-simplex.}
\label{fig:TwoSimplices}
\end{figure}

Adding noise -- random bit flips with probability $\alpha$ -- to the data separates the effective probability of ``heads'',
\begin{equation} \label{eqQofP}
q = \alpha + p(1-2\alpha),
\end{equation}
from the true probability $p$.  Quite a lot of the complications that ensue can be understood as stemming from a single underlying schizophrenia in the problem:  there are now \emph{two} relevant probability simplices, one for $q$ and one for $p$ (see Figure \ref{fig:TwoSimplices}).

One part of the problem (the data, and therefore the likelihood function) essentially live on the $q$-simplex.  The other parts (the parameter to be estimated, and therefore the risk function) live on the $p$-simplex.  The simplest, most obvious estimator is linear inversion,
$$\hat{q} = \frac{n}{N},$$
which implies (by inversion of Eq. \ref{eqQofP}),
\begin{equation}\label{equation Robin wanted labeled :p}
\hat{p} = \frac{\hat{q}-\alpha}{1-2\alpha}.
\end{equation}
But if $n<\alpha N$ or $n>\left(1-\alpha\right)N$, then $\hat{p}$ is negative or greater than 1.  This is patently absurd.  It is also clearly suboptimal, as there is no advantage to assigning an estimate that lies outside the (convex) set of valid states.  Finally, it guarantees infinite expected risk -- which is to say that it is quantitatively \emph{very} suboptimal.

Maximum likelihood does not improve matters much.  A bit of algebra shows that, just as for the noiseless coin, ML coincides with linear inversion -- \emph{if} $\pML$ is a valid probability.  Otherwise, $\cL(p)$ achieves its maximum value on the boundary of its domain, at 0 or 1.
\begin{equation*}
\pML = \left\{\begin{array}{cc} 0 & \mathrm{if\ }n<\alpha N \\
  \frac{n-\alpha N}{N(1-2\alpha)} & \mathrm{if\ }\alpha N \leq n \leq N(1-\alpha) \\
  1 & \mathrm{if\ }n>N(1-\alpha)\end{array}\right.
\end{equation*}
The ML estimator is always a valid probability (by construction, since the domain of $\cL(p)$ is the simplex).  However, like linear inversion, it is still clearly suboptimal.  It is never risk-minimizing to assign $\hat{p}=0$ unless we are \emph{certain} that $p$ truly is zero.  Moreover, the expected entropy risk is still infinite under all circumstances, since $\hat{p}=0$ occurs with nonzero probability and $R(0;p)=\infty$.

For the noiseless coin the \emph{only} way that $\pML=0$ is if $n=0$ is observed.  For the noisy coin, $\pML=0$ for a whole range of data -- and with probability close to 50\% when $p\ll\frac{1}{\sqrt{N}}$.  In the noiseless context, adding fictious observations solved this problem quite well, generating ``add $\beta$'' estimators with near-optimal performance.  Unfortunately, the addition of fictitious observations
$$n_k \to n_k + \beta$$
doesn't have the same effect for the noisy coin.  If we have $n < \alpha N - \beta$ (which is possible and even probable), then adding $\beta$ has no effect.  Linear inversion still gives $\hat{p}<0$, and ML still gives $\hat{p}=0$.

The core problem here is sampling mismatch.  We sample from $q$, but the risk is determined by $p$.  ML takes no account of the risk function, and neither does our attempt to hedge the estimate by adding fake samples.  Both are entirely $q$-centric.  A mechanism that acts directly on the $p$-simplex is needed, to force $\hat{p}$ away from the dangerous boundaries ($\hat{p}=0$ and $\hat{p}=1$, where the risk diverges).

\begin{figure}[th]
\begin{tabular}{l}
\includegraphics[width=3in]{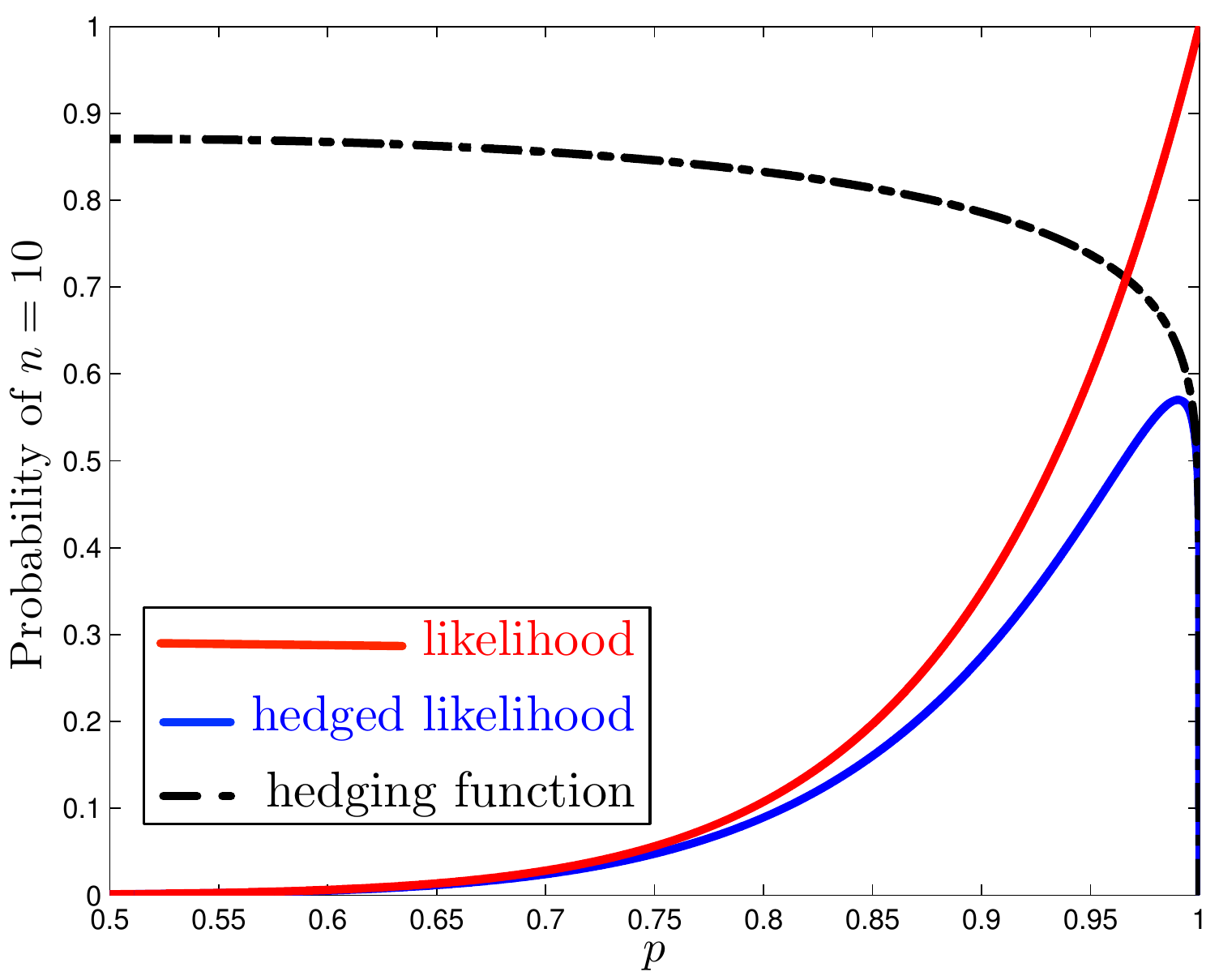} \\ \includegraphics[width=3in]{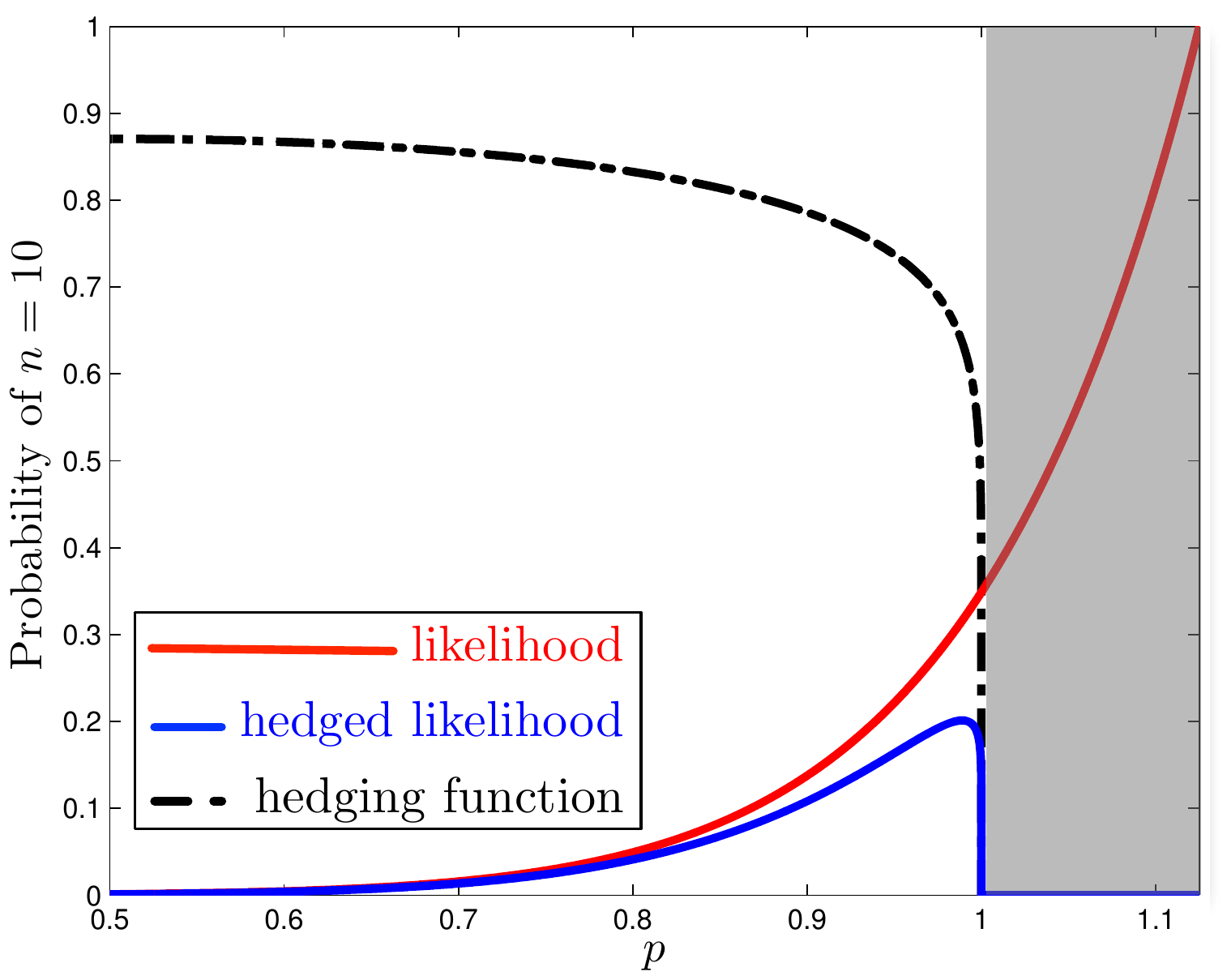}
\end{tabular}
\caption{ Illustration of a hedging function ($\beta=0.1$), and its
effect on the likelihood function for noiseless (top) and noisy (bottom)
coins.  In the top plot, we show the hedging function over the 2-simplex
(dotted black line), the likelihood function for an extreme data set
comprising 10 heads and 0 tails (red line), and the corresponding
\emph{hedged} likelihood (blue line).  The lower plot shows the same
functions for a noisy coin with $\alpha=0.1$.  The shaded regions are
outside the $p$-simplex (and therefore forbidden), but correspond to
valid $q$ values.  Note that the \emph{unconstrained} maximum of $\cL$
lies in the forbidden region where $p>1$, and therefore the maximum of
the constrained likelihood is on the boundary ($p=1$) -- a pathology
that hedging remedies.}
\label{fig:hedgingfunction}
\end{figure}

One simple way to do this is to modify $\cL(p)$, multiplying it by a ``hedging function'' \cite{BlumeKohout2010Hedged} $h(p) = \prod_k{p_k^\beta} = p^\beta(1-p)^\beta$,
$$\cL(p) \to \cL'(p) = p^\beta(1-p)^\beta\cL(p),$$
and define $\hat{p}_\beta = \mathrm{argmax}\cL'(p)$ (see Fig. \ref{fig:hedgingfunction}).  For a noiseless coin, this is identical to adding $\beta$ fictitious observations of each possible event -- but for $\alpha>0$, they are not equivalent.  The hedging function modification is sensitive to the $p_k=0$ boundary of the simplex, and inexorably forces the maximum of $\cL'(p)$ away from it (since $\cL'(p)$ remains log-convex, but equals zero at the boundary).

The HML estimator is given by
$$\hat{p}_\beta = \frac{\hat{q}_\beta-\alpha}{1-2\alpha},$$
where $\hat{q}_\beta$ is the zero of the cubic polynomial
\begin{equation*} \label{eq:NoisyHMLECubic}
(N+2\beta)\hat q^3-(N+n+3\beta)\hat q^2+(n+\beta+N\alpha-N\alpha^2)\hat q+n\alpha^2-n\alpha.
\end{equation*}
that lies in $[\alpha,1-\alpha]$.

\begin{figure}[th]
\includegraphics[width=3in]{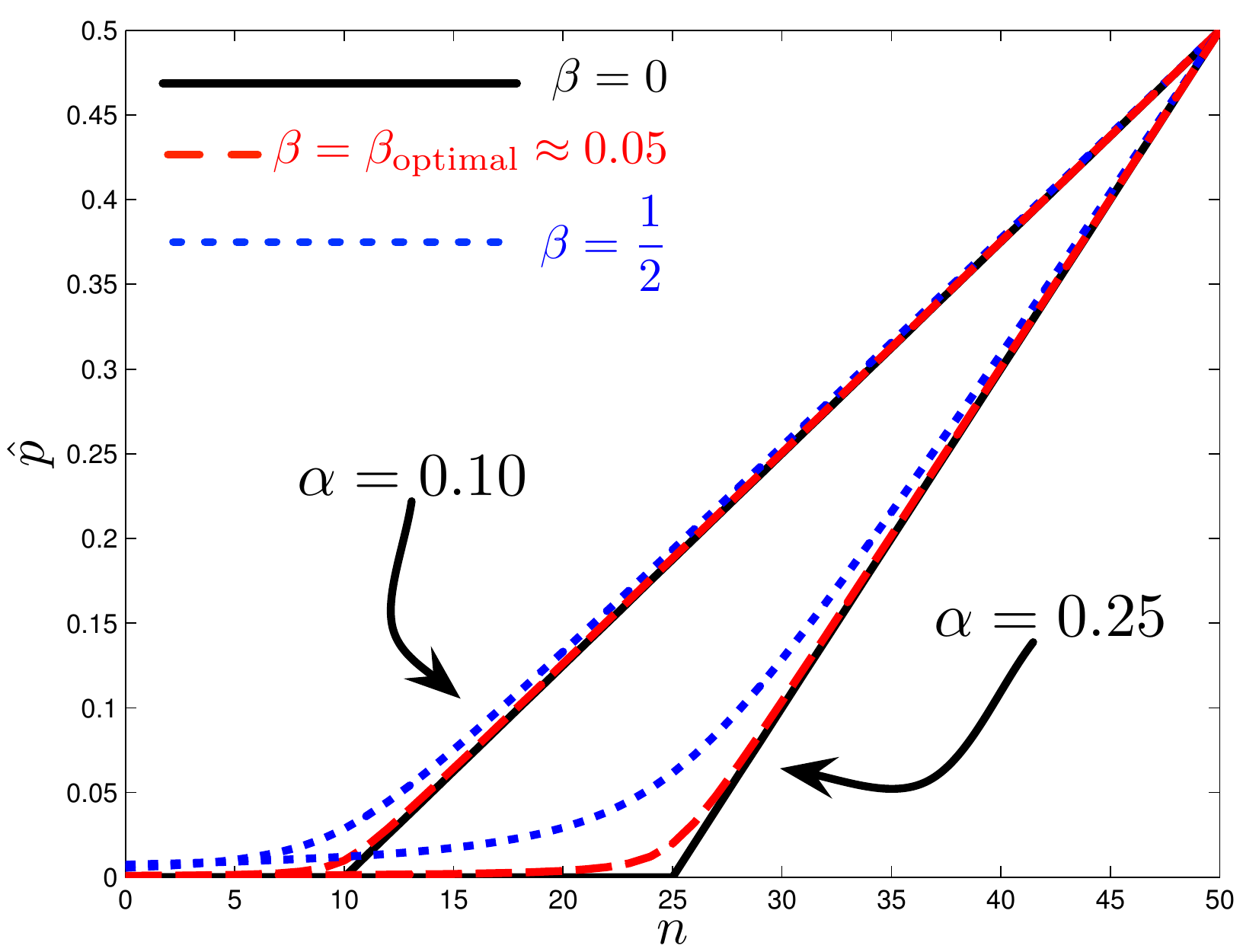}
\caption{HML (hedged maximum likelihood) estimators $\hat{p}_\beta(n)$ are shown for several values of $\beta$, and compared with the maximum likelihood (ML) estimator $\pML(n)$.  $N=100$ in all cases.  Whereas the ML estimator is linear in $n$ until it encounters $p=0$, the hedged estimator smoothly approaches $p=0$ as the data become more extreme.  Increasing $\beta$ pushes $\hat{p}$ away from $\hat{p}=0$.}
\label{fig:HMLEstimators}
\end{figure}

Figure \ref{fig:HMLEstimators} illustrates how $\hat{p}_\beta$ and $\pML$ depend on $\frac{n}{N}$.  When $\pML$ is far from the simplex boundary (0 and 1), hedging has relatively little effect.  In fact, hedging yields an approximately linear estimator akin to ``add $\beta$''.  But as $\pML$ approaches 0 or 1, the effect of hedging increases.  When $\pML$ intersects 0, at $n=\alpha N$, $\hat{p}_{\beta} = O(1/\sqrt{N})$.  This fairly dramatic shift occurs because the likelihood function is approximately Gaussian, with a maximum at $p=0$ and a width of $O(1/\sqrt{N})$.  $\cL(p)$ declines rather slowly from $p=0$, and $p=O(1/\sqrt{N})$ is not \emph{substantially} less likely than $p=0$, so the hedging imperative to avoid $\hat{p}=0$ pushes the maximum of $\cL'(p)$ far inside the simplex.

\begin{figure}[th]
\begin{tabular}{c}
\includegraphics[width=3in]{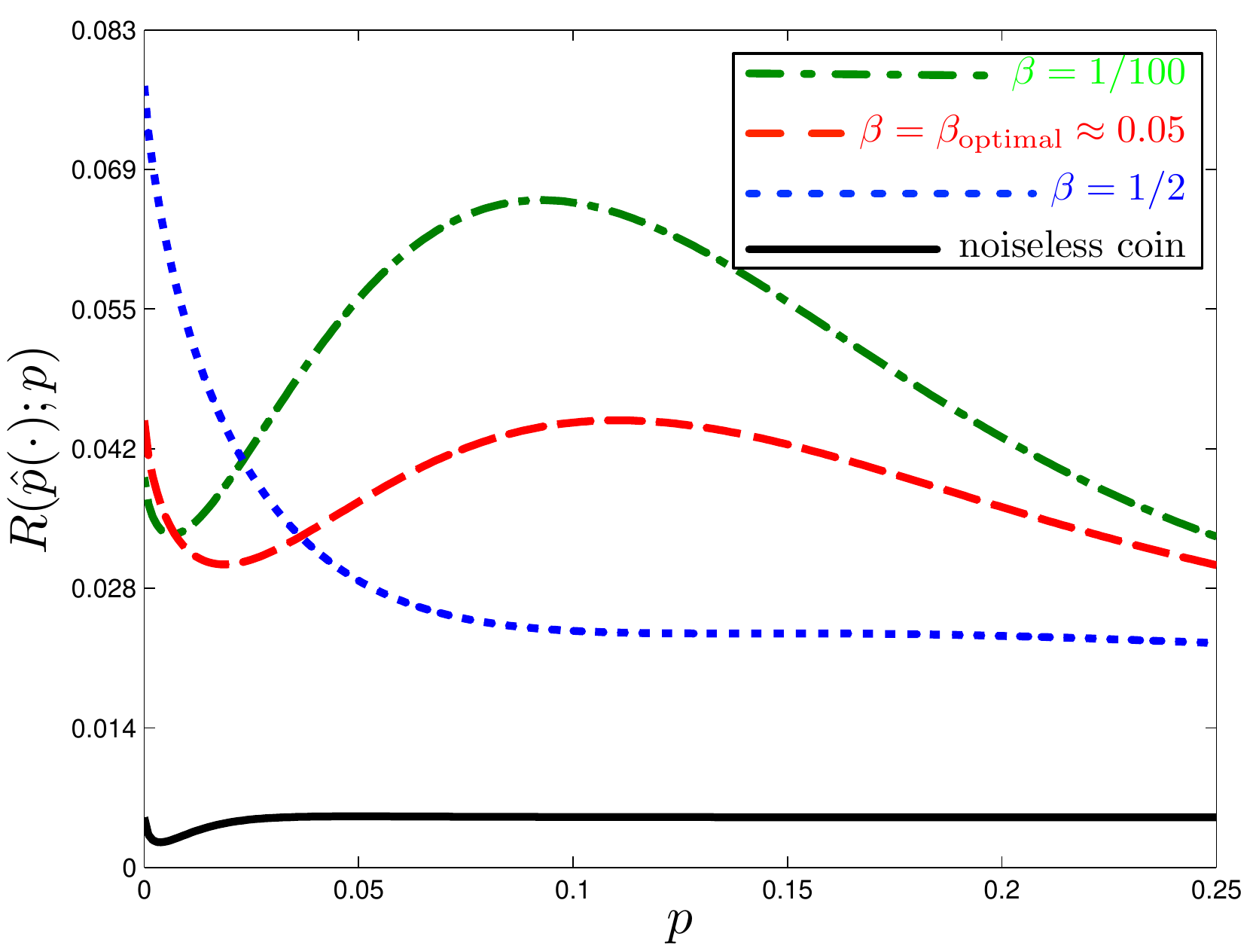} \\ \includegraphics[width=3in]{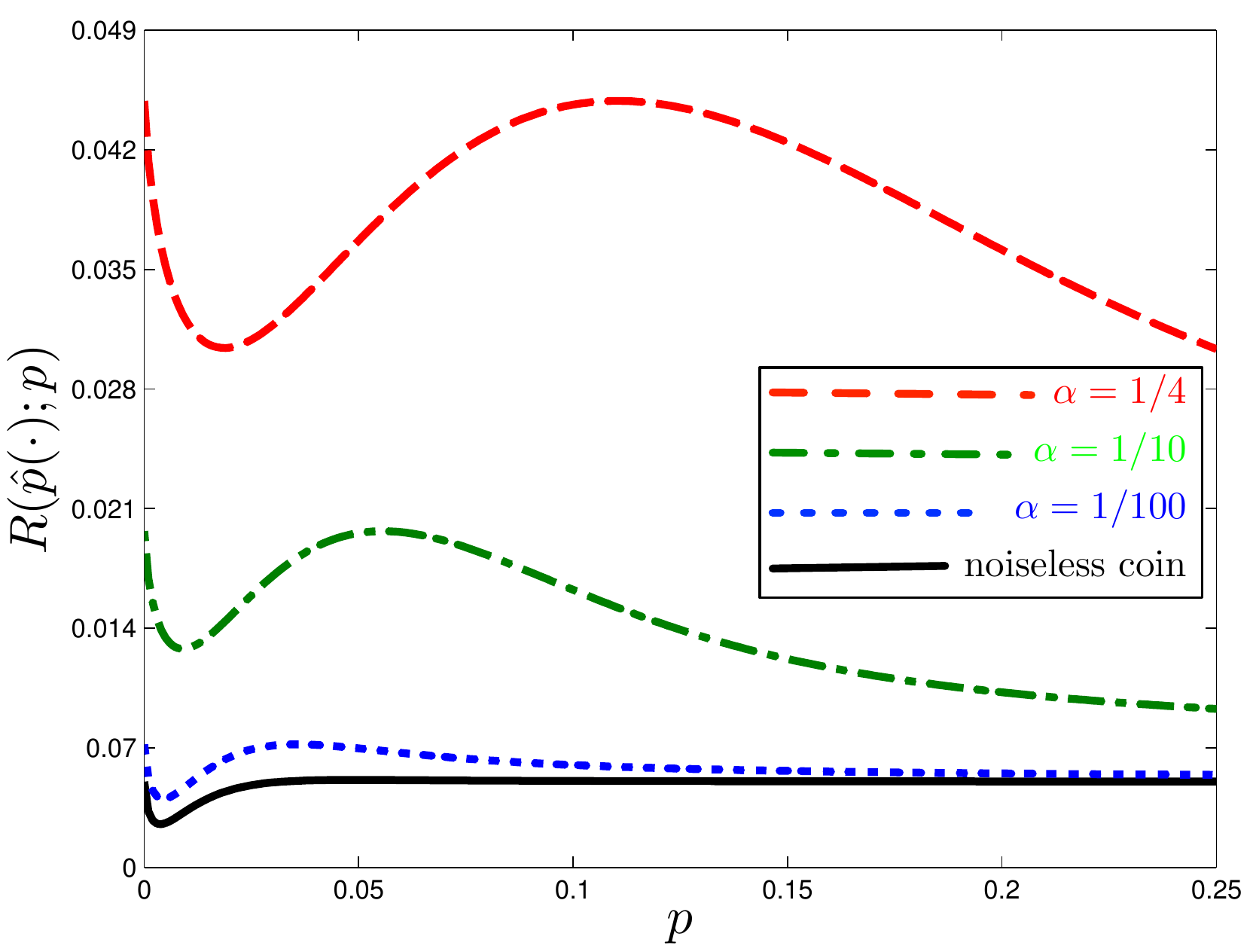}
\end{tabular}
\caption{The risk profile $R(p)$ is shown for several HML estimators.  $N=100$ in all cases.  In the top plot, a noisy coin with $\alpha=1/4$ has been estimated using three different HML estimators.  The optimal $\beta$ balances boundary risk against interior risk.  Increasing $\beta$ increases boundary risk, while decreasing it increases interior risk.  The bottom plot shows the risk of optimal HML estimators for $\alpha=0,1/100,1/10,1/4$.  Risk approaches $O(1/\sqrt{N})$ for noisy coins, vs. $O(1/N)$ for noiseless coins.}
\label{fig:HMLRisk}
\end{figure}

How accurate are these hedged estimators?  Figure \ref{fig:HMLRisk} shows the \emph{pointwise} average risk as a function of the true $p$, for different amounts of noise ($\alpha$) and hedging ($\beta$).  For the noiseless ($\alpha=0$) coin, $\beta=1/2$ yields a nearly flat risk profile given by $R(p) \approx 1/2N$.  In contrast, hedged estimators for the noiseless coin yield similar profiles that rise from $O(1/N)$ in the interior to a peak of $O(1/\sqrt{N})$ around $p = O(1/\sqrt{N})$.  Risk \emph{at} the $p=0$ boundary depends on $\beta$, and may be either higher or lower than the peak at $p = O(1/\sqrt{N})$.

At first glance, this behavior suggests a serious flaw in the hedged estimators.  The peak around $p \approx O(1/\sqrt{N})$ is of particular concern, since in all cases the risk is $O(1/\sqrt{N})$ there.  But in fact, this behavior is generic for the noisy coin.  Minimax estimators have
similar $O(1/\sqrt N)$ errors, and hedged estimators turn out to perform quite well.  However, they are \emph{not} minimax, or even close to it!  As we shall show in the next section, the noisy coin's ``intrinsic risk'' profile is far from flat.  The minimax estimator attempts to flatten it -- at substantial cost.

\section{Minimax estimators for the noisy coin}

Some simple estimators (such as ``add $1/2$'') are nearly minimax for the noiseless coin.  This is not true for the noisy coin in general, because (as we shall see) the minimax estimators are somewhat pathological.  So we used numerics to find good approximations to minimax estimators for noisy coins.

Minimax-Bayes duality permits us to search over priors rather than estimators.  Each prior $\pi(p)\diff p$ defines a Bayesian mean estimator $\hat{p}_\pi(n)$, which is Bayes for $\pi(p)$.  Its risk profile $R_{\hat{p}_\pi}(p)$ provides both upper and lower bounds on the minimax risk (Eq. \ref{eq:RiskBounds}).

As is often the case for discretely distributed data, the minimax priors for coins appear to always be discrete \cite{Lehmann1998Theory}.  We searched for least favorable priors (holding $N$ and $\alpha$ fixed) using the algorithm of Kempthorne \cite{Kempthorne1987Numerical}.  We defined a prior with a few support points, and let the location and weight of the support points vary in order to maximize the Bayes risk.  Once the optimization equilibrated, we added new support points at local maxima of the risk, and repeated this process until the algorithm found priors for which the maximum and Bayes risk coincided to within $10^{-6}$ relative error.

\begin{figure}[th]
%\begin{tabular}{ll}
\includegraphics[width=3in]{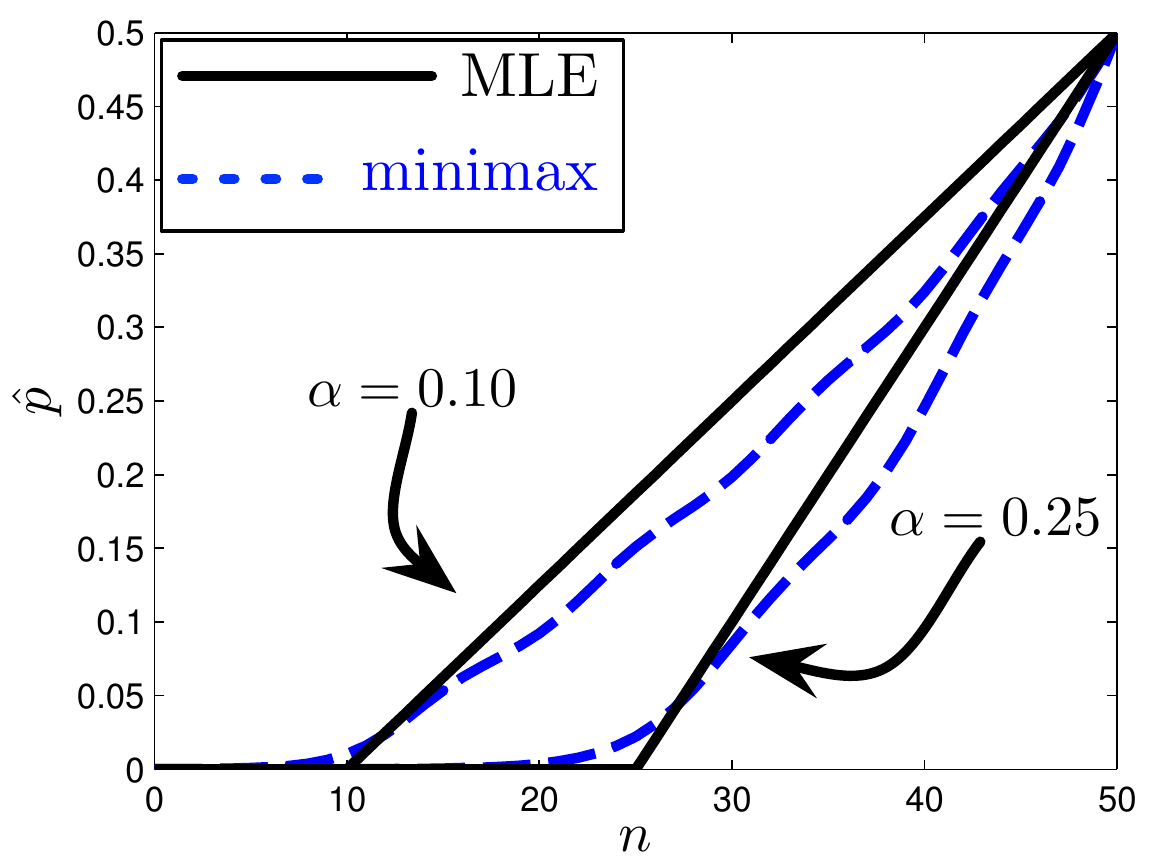}
%\end{tabular}
\caption{Minimax and ML estimators are shown for $N=100$ and $\alpha=1/10,1/4$.  Note that the minimax estimator is grossly biased in the interior -- a pathological result of the mandate to minimize \emph{maximum} risk at all costs.}
\label{fig:MinimaxEstimators}
\end{figure}

\begin{figure}[th]
\begin{tabular}{lc}
\includegraphics[width=3in]{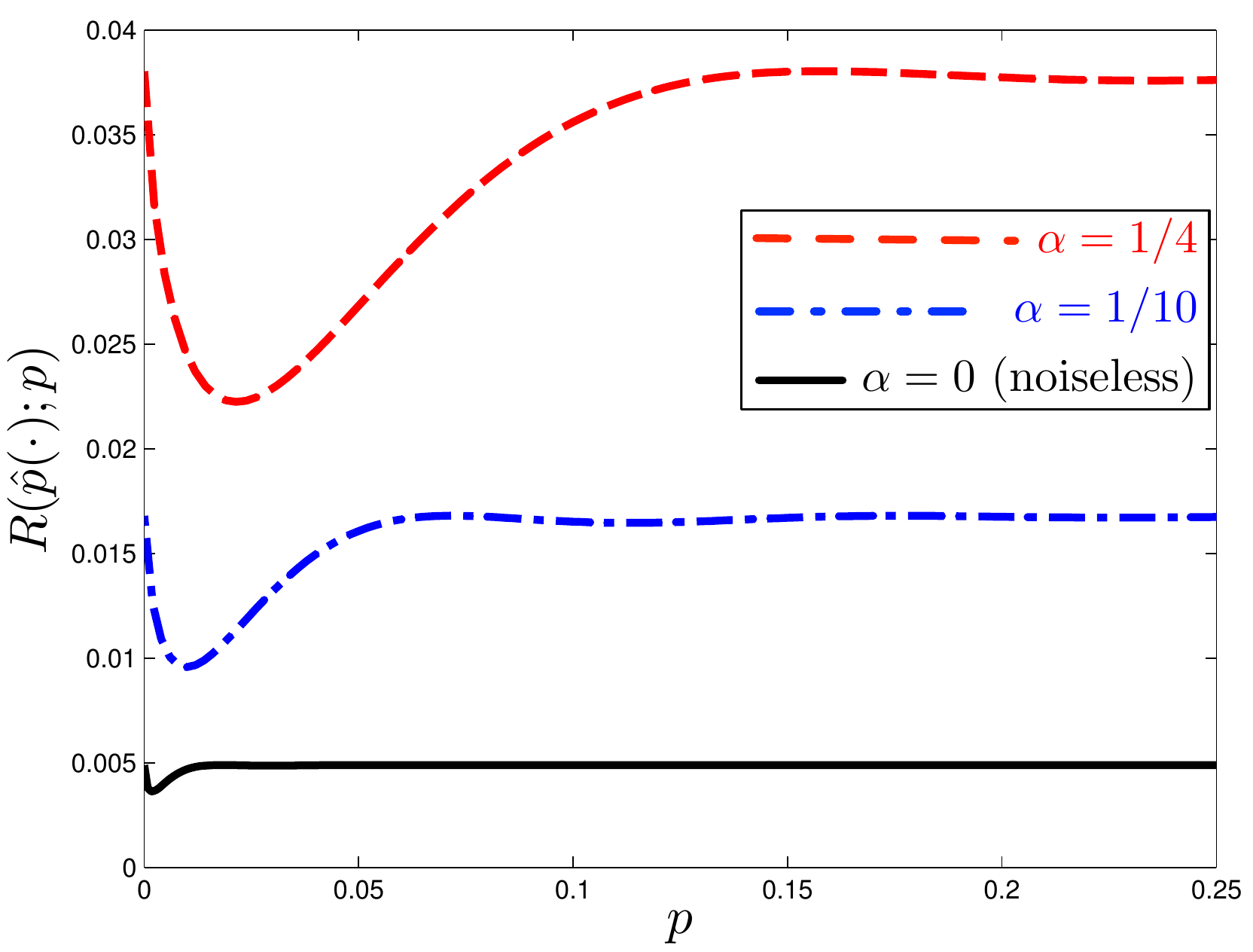} \\ \includegraphics[width=3in]{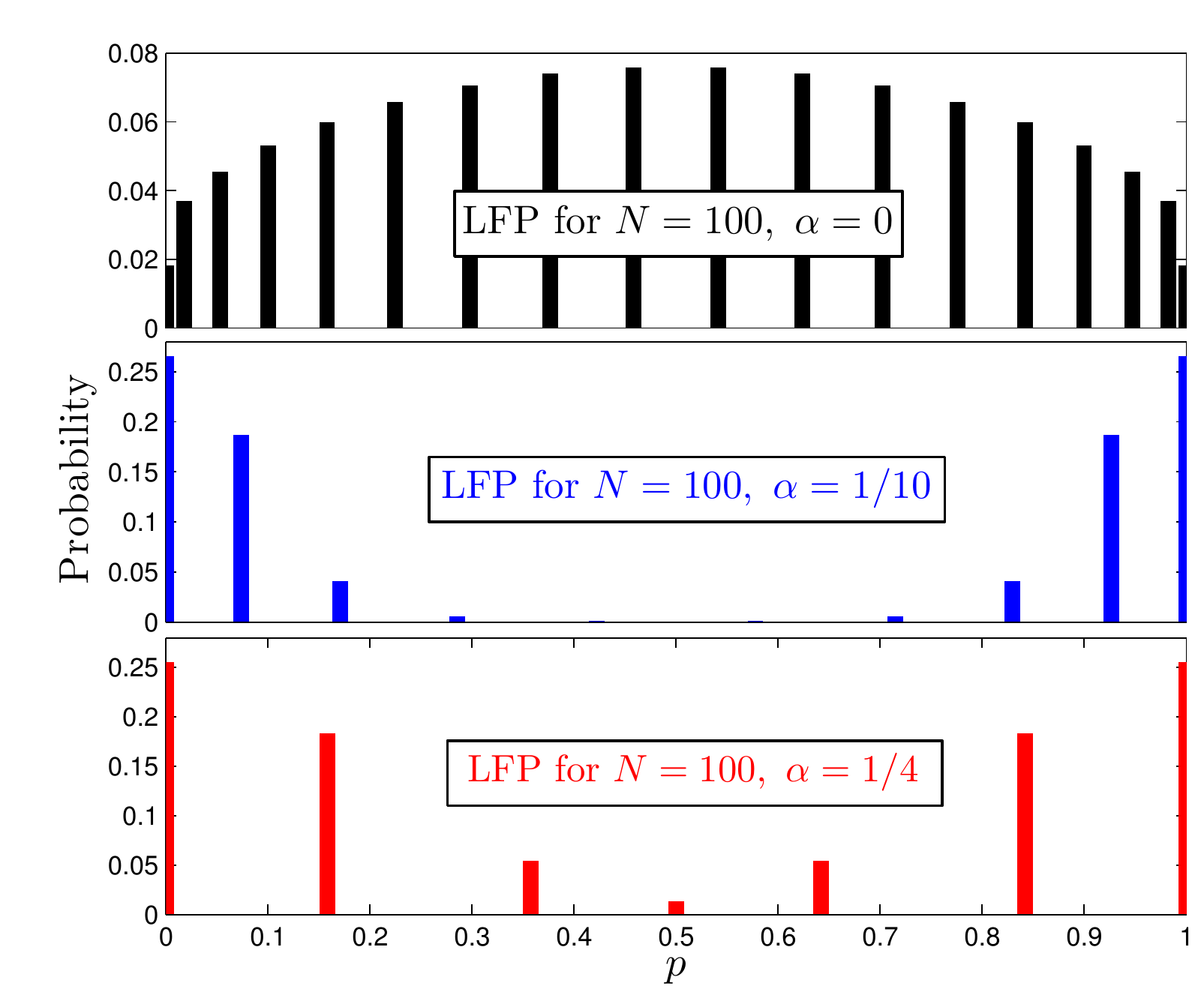}
\end{tabular}
\caption{Above (top), the risk profile $R(p)$ is shown for the minimax estimators of Figure \ref{fig:MinimaxEstimators} ($N=100$; $\alpha=1/10,1/4$) and for a noiseless coin (also $N=100$).  No estimator can achieve lower risk across the board -- but the minimax estimator's risk is very high in the interior ($p>1/\sqrt{N}$) compared with HML estimators.  The $O(1/\sqrt{N})$ risk of HML estimators is intrinsic to noisy coins.  Below, we show the least favorable priors whose Bayes estimators have the risk profiles above.  The support points of the discrete priors occur, as they must, where the risk achieves its maximum value.  Thus, the Bayes risk is equal to maximum risk (within a relative tolerance of $10^{-6}$).}
\label{fig:MinimaxRisk}
\end{figure}

%\begin{figure}[th]
%\begin{tabular}{ll}
%\includegraphics[width=3in]{Fig_2DSimplexMismatch} & \includegraphics[width=2in]{Fig_3DSimplexMismatch}
%\end{tabular}
%\caption{\texttt{This is a sample LFP.  With parameters chosen basically on the grounds that we want the most nice-looking LFP.}}
%\label{fig:LFP}
%\end{figure}

Figure \ref{fig:MinimaxEstimators} illustrates minimax estimators for several $N$ and $\alpha$, while Figure \ref{fig:MinimaxRisk} shows the resulting risk profiles.  The minimax risk is $O(1/\sqrt{N})$ -- \emph{not} $O(1/N)$ as for the noiseless coin.  The risk is clearly dominated by points near the boundary, and we find that the LFPs typically place almost all their weight on support points within a distance $O(1/\sqrt{N})$ of the boundary ($p=0$ and $p=1$).  As a result of this severe weighting toward the boundary, the minimax estimators are highly biased toward $p \approx 1/\sqrt{N}$ -- not just when $p$ is close to the boundary (when bias is inevitable) but also when $p$ is in the interior!  This effect is truly pathological, although it can easily be explained.  Low risk, of order $1/N$, can easily be achieved in the interior.  However, the minimax estimator seeks at all costs to reduce the \emph{maximum} risk, which is achieved near $p\approx 1/\sqrt{N}$.  By biasing heavily toward $p\approx1/\sqrt{N}$, the estimator achieves slightly lower maximum risk\ldots at the cost of dramatically increasing its interior risk from $O(1/N)$ to $O(1/\sqrt{N})$.

The preceding analysis made use of an intuitive notion of pointwise ``intrinsic risk'' -- i.e., a lower bound $R_{\mathrm{min}}(p)$ on the expected risk for any given $p$.  Formally, no such lower bound exists.  We can achieve $R(p')=0$ for \emph{any} $p'$, simply by using the estimator $\hat{p}=p'$.  But we can rigorously define something very similar, which we call \emph{bimodal risk}.

The reason that it's not practical to achieve $R(p')=0$ at any given $p'$ is, of course, that $p'$ is unknown.  We must take into account the possibility that $p$ takes some other value.  Least favorable priors are intended to quantify the risk that ensues, but a LFP is a property of the entire problem, not of any particular $p'$.  In order to quantify ``how hard is a particular $p'$ to estimate,'' we consider the set of bimodal priors,
$$\pi_{w,p',p''}(p) = w\delta(p-p') + (1-w)\delta(p-p''),$$
and maximize Bayes risk over them.  We define the bimodal risk of $p'$ as
$$R_{2}(p') = \max_{w,p''}{\Ravg_{\pi_{w,p',p''}}}.$$
The bimodal risk quantifies the difficulty of distinguishing $p'$ from \emph{just one} other state $p''$.  As such, it is always a lower bound on the minimax risk.

\begin{figure}[th]
%\begin{tabular}{ll}
\includegraphics[width=3in]{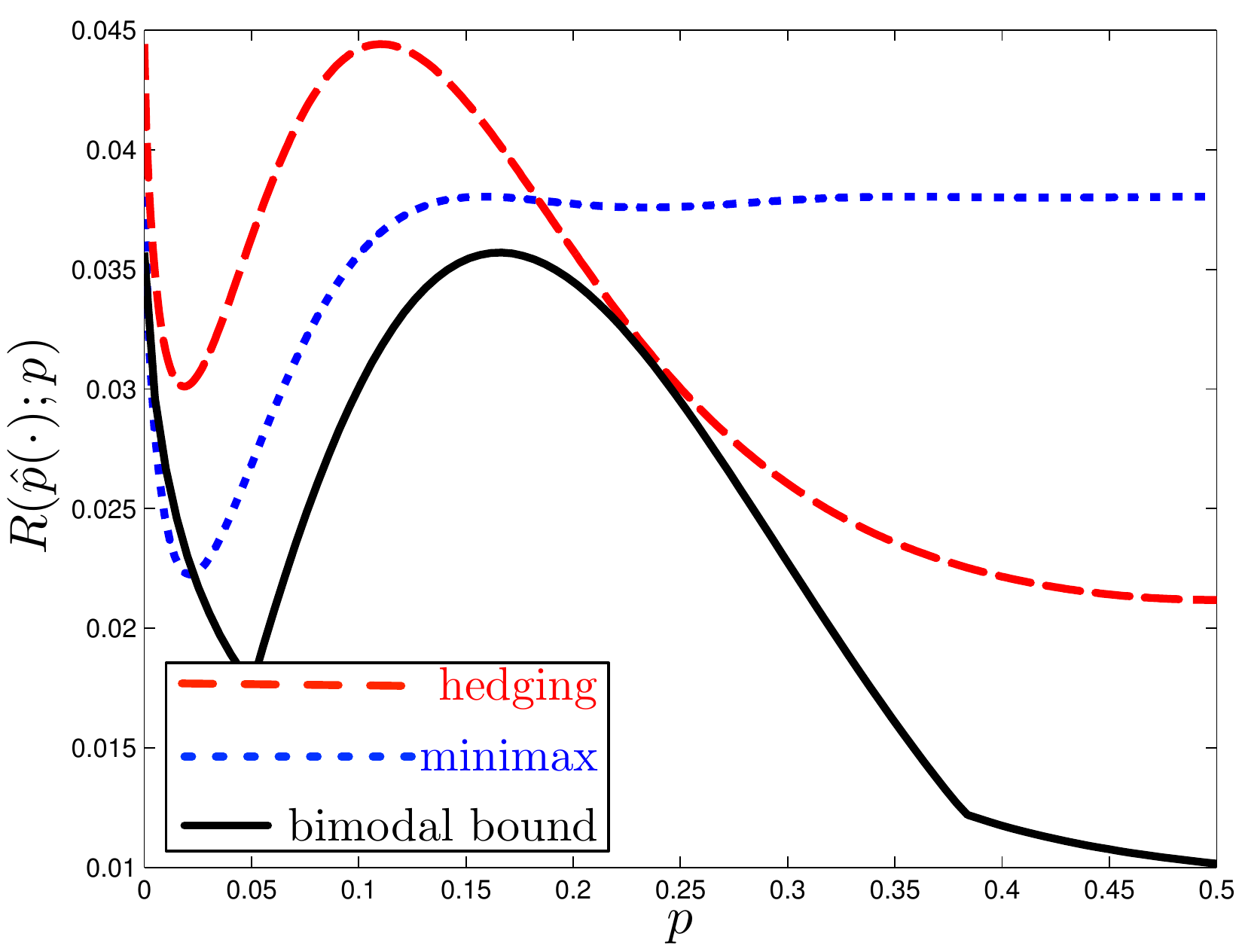}
%\end{tabular}
\caption{The risk profiles of the optimal HML estimator (red) and the minimax estimator (blue) are compared with the bimodal lower bound $R_2(p)$.  Note that while the HML maximum risk exceeds the minimax risk (as it must!), it is competitive -- and HML is \emph{far} more accurate in the interior.  The bimodal lower bound supports the conjecture that HML is a good compromise, since the HML risk exceeds the bimodal bound by a nearly constant factor.}
\label{fig:OptimalRisks}
\end{figure}

Figure \ref{fig:OptimalRisks} compares the bimodal risk to the pointwise risk achieved by the minimax and [optimal] HML estimators.  Note that the bimodal risk function is a strict lower bound (at every point) for the minimax risk, but not for the pointwise risk of any estimator (including the
minimax estimator).  However, every estimator exceeds the bimodal risk at at least one point, and almost certainly at \emph{many} points.  Figure \ref{fig:OptimalRisks} confirms that the noisy coin's risk is dominated by the difficulty of distinguishing $p\approx 1/\sqrt{N}$ from $p=0$.  States deep inside the simplex are far easier to estimate, with an expected risk of $O(1/N)$.

A simple analytic explanation for this behavior can be obtained by series expansion of the Kullback-Leibler divergence,
$$KL(p+\epsilon;p) \approx \frac{\epsilon^2}{2p(1-p)}.$$
The typical error, $\epsilon$, in the estimator is $O(1/\sqrt{N})$, so as long as both $p$ and $1-p$ are bounded away from zero, the typical entropy risk is $O(1/N)$.  But as $p$ approaches 0 (or 1), this approximation diverges.  If we consider $p=1/\sqrt{N}$, then we get
$$KL(p+\epsilon;p) \approx \frac{\epsilon^2}{2p} = O(1/\sqrt{N}).$$
At $p=0$, the series expansion fails, but a one-sided approximation (with $\epsilon$ strictly greater than 0) gives the same result,
$$KL(\epsilon;0) \approx \epsilon = O(1/\sqrt{N}).$$
For the noiseless coin, this does not occur because the typical error is not always $O(1/\sqrt{N})$.  Instead, it scales with the inverse of the Fisher information,
$$\epsilon \approx \sqrt{\frac{2p(1-p)}{N}},$$
and the $p$-dependent factor neatly cancels its counterpart in the Kullback-Leibler divergence, which produces the nice flat risk profile seen for the noiseless coin.  The underlying problem for the noisy coin is -- again -- sampling mismatch.  Because we are observing $q$ and predicting $p$, the two factors do not cancel out.

\section{Good estimators for the noisy coin}

Minimax is an elegant concept, but for the noisy coin it does not yield ``good'' estimators.  In a single-minded quest to minimize the maximum risk, it yields wildly biased estimates in the interior of the simplex.  This is reasonable only in the [implausible] case where $p$ is truly selected by an adversary (in which case the adversary would almost always choose $p$ near the boundary, out of sheer bloody-mindedness).  In the real world, robustness against adversarial selection of $p$ is good, but should not be taken to absurd limits.

This leaves us in need of a quantitative criterion for ``good'' estimators.  Ideally, we would like an estimator that achieves (or comes close to achieving) the ``intrinsic'' risk for every $p$.  The bimodal risk $R_2(p)$ provides a reasonably good proxy -- or, more precisely, a lower bound -- for intrinsic risk (in the absence of a rigorous definition).  This is not a precise quantitative framework, but it does provide a reasonably straightforward criterion:  we are looking for an estimator that closely approaches the bimodal risk profile.

\begin{figure}[t]
%\begin{tabular}{ll}
\includegraphics[width=3in]{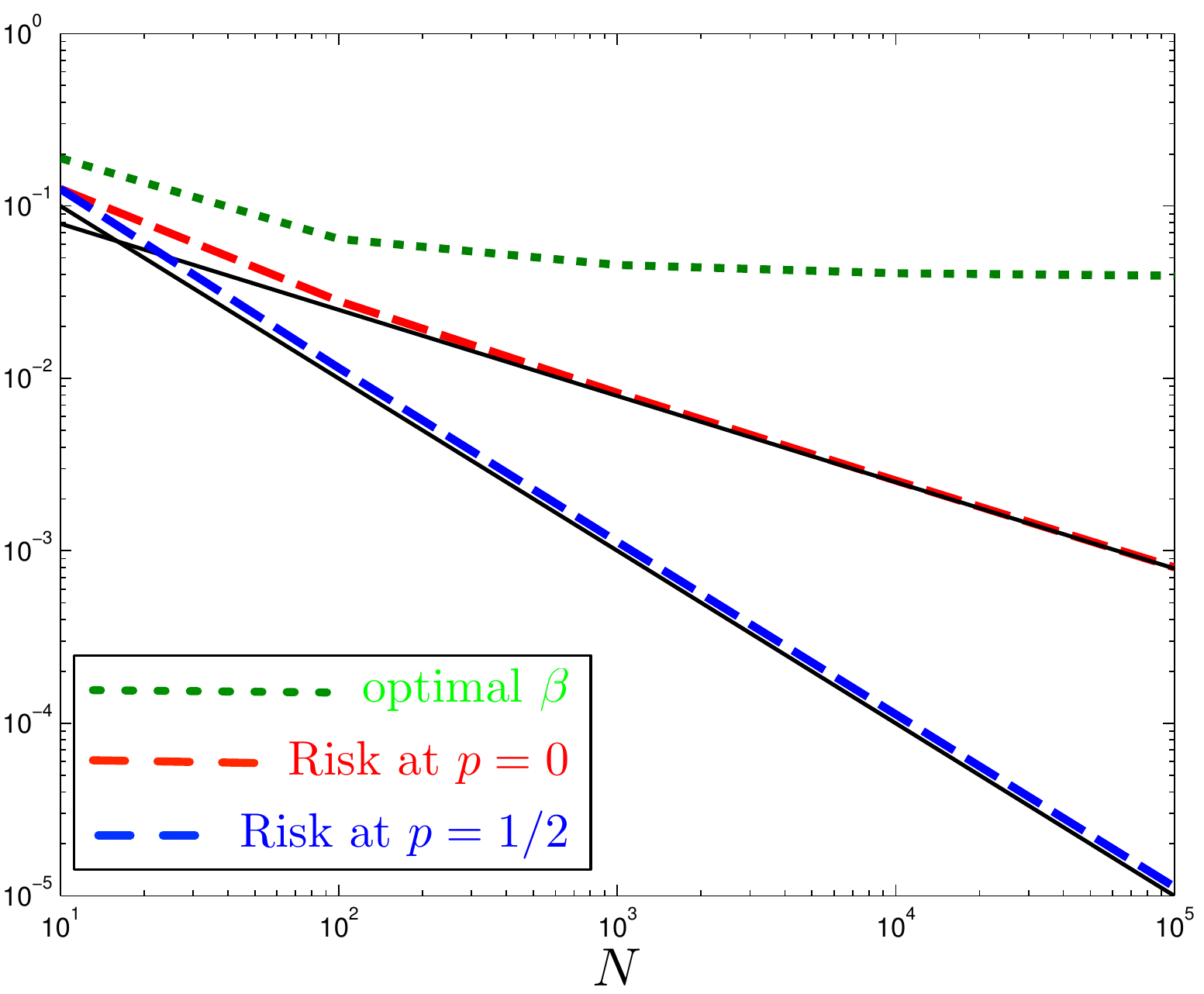}
%\end{tabular}
\caption{Here, we show how three aspects of optimal HML estimation behave as $N$ increases toward infinity.  They are: (1) the risk deep in the bulk (blue); (2) the bimodal risk at the boundary (red); and (3) the optimal value of $\beta$ (green).  The data shown here, $N=10\ldots10^5$, correspond to a noise level of $\alpha=1/100$, but other values of $\alpha$ produce qualitatively identical results.  The ``optimal'' $\beta$ is the one that achieves the smallest maximum risk, by balancing the risk at $p=0$ against its maximum value in the interior of the simplex.  Its value decreases with $N$ for small $N$, and for $N \gg \alpha^{-1}$ it asymptotes to $\beta_{\mathrm{optimal}}\approx 0.0389$ (see also Fig. \ref{fig:OptimalBeta2}).  The HML estimator's expected risk varies quite dramatically with $p$ (as expected).  At $p=0$, the risk asymptotes to $R(0) \approx 1/4\sqrt{N}$.  At $p=1/2$, it asymptotes to $R(1/2) \approx 1/N$.  The $O(1/\sqrt{N})$ expected risk is seen only for $p \in [0,O(1)/\sqrt{N}]$.}
\label{fig:OptimalBeta}
\end{figure}

\begin{figure}[th]
%\begin{tabular}{ll}
\includegraphics[width=3in]{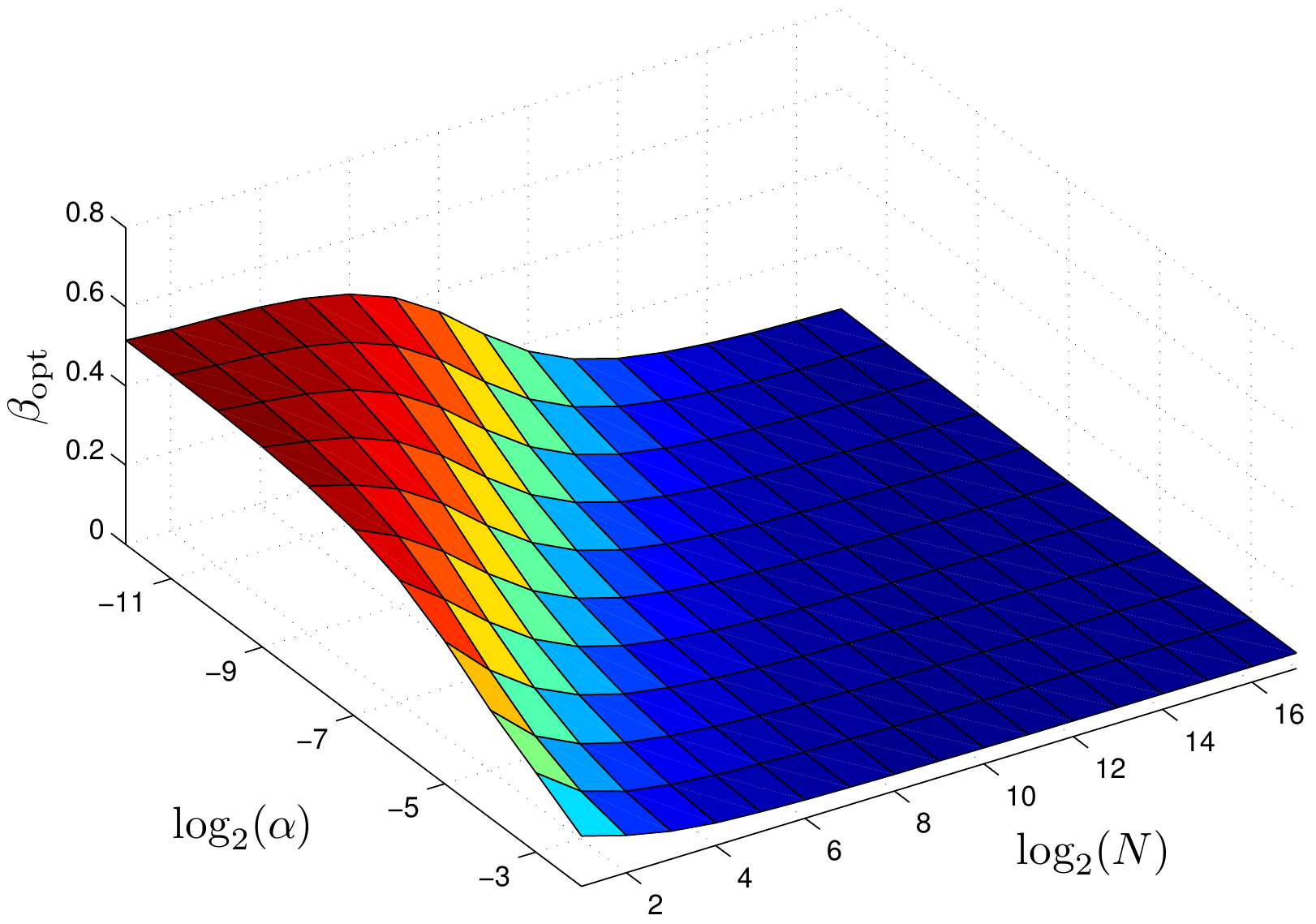}
%\end{tabular}
\caption{The optimal value of $\beta$ for $N=2\dots2^{17}$ and $\alpha=2^{-12}\ldots2^{-2}$.  It approaches the optimal noiseless value $\approx 1/2$ when $N\ll\alpha^{-1}$.  It rapidly declines and at roughly $N\approx\alpha^{-1}$ it is well within $10^{-2}$ of what appears to be its asymptotic value $\beta_{\mathrm{optimal}}\approx 0.0389$.}
\label{fig:OptimalBeta2}
\end{figure}

Hedged estimators are a natural ansatz, but we need to specify $\beta$.  Whereas the noiseless coin is fairly accurately estimated by $\beta=1/2$ for all $N$, the optimal value of $\beta$ varies with $N$ for noisy coins.  Local maxima of the risk are located at $p=0$ and at $p\approx1\sqrt{N}$, one or both of which is always the global maximum.  So, to choose $\beta$, we minimize maximum risk by setting them equal to each other.  Figure \ref{fig:OptimalBeta} shows the optimal $\beta$ as a function of $N$, for a representative value of $\alpha$, which approaches $\beta_{\mathrm{optimal}} \approx 0.0389$ for large $N$.  This value is obtained for a large range of $\alpha$'s, as shown in Figure \ref{fig:OptimalBeta2}.  Finally, Figure \ref{fig:OptimalRisks} compares the risk profile of the optimal hedging estimator with (i) bimodal risk, and (ii) the risk of the minimax estimator.  We conclude that while optimal hedging estimators probably do not offer strictly optimal performance, they are (i) easy to specify and calculate, (ii) far better than minimax estimators for almost all values of $p$, and (iii) relatively close to the lower bound defined by bimodal risk.

\section{Applications and discussion}

Our interest in noisy coins (and dice) stems largely from their wide range of applications.  Our original motivation came from quantum state and process estimation (see, e.g. \cite{BlumeKohout2010Hedged}, as well as reference therein).  However, noisy coins appear in several other important scientific contexts.

\subsection{Randomized response}

Social scientists often want to know the prevalence of an embarrassing or illegal habit -- e.g., tax evasion, abortion, infidelity -- in a population.  Direct survey yields negatively biased results, since respondents often do not trust guarantees of anonymity.  Randomized response \cite{Chaudhuri1988Randomized} avoids this problem by asking each respondent to flip a coin of known bias $\alpha$ and invert their answer if it comes up heads.  An adulterer or tax cheat can answer honestly, but claim (if confronted) that they were lying in accordance with instructions.  The scientist is left with the problem of inferring the \emph{true} prevalence ($p$) from noisy data -- precisely the problem we have considered here.  Interval estimation (e.g., confidence intervals) are a good alternative, but if the study is intended to yield a point estimate, our analysis is directly applicable.

\subsection{Particle detection}

Particle detectors are ubiquitous in physics, ranging from billion-dollar detectors that count neutrinos and search for the Higgs boson, to single-photon photodetectors used throughout quantum optics.  The number of counts is usually a Poisson variable, which differs only slightly from the binomial studied here.  \emph{Background counts} are an unavoidable issue, and lead to an estimation problem essentially identical to the noisy coin.  Particle physicists have argued extensively (see, e.g., \cite{MandelkernSS02,RoePRD99}, and references therein) over how to proceed when the observed counts are \emph{less} than the expected background.  Our proposed solution (if a point estimator is desired -- region estimation, as in \cite{FeldmanCousins}, may be a better choice) is to use HML or another estimator with similar risk profile -- and to be aware that the existence of background counts has a dramatic effect on the expected risk.

\subsection{Quantum state (and process) estimation}

\begin{figure}[t]
\includegraphics[width=2.5in]{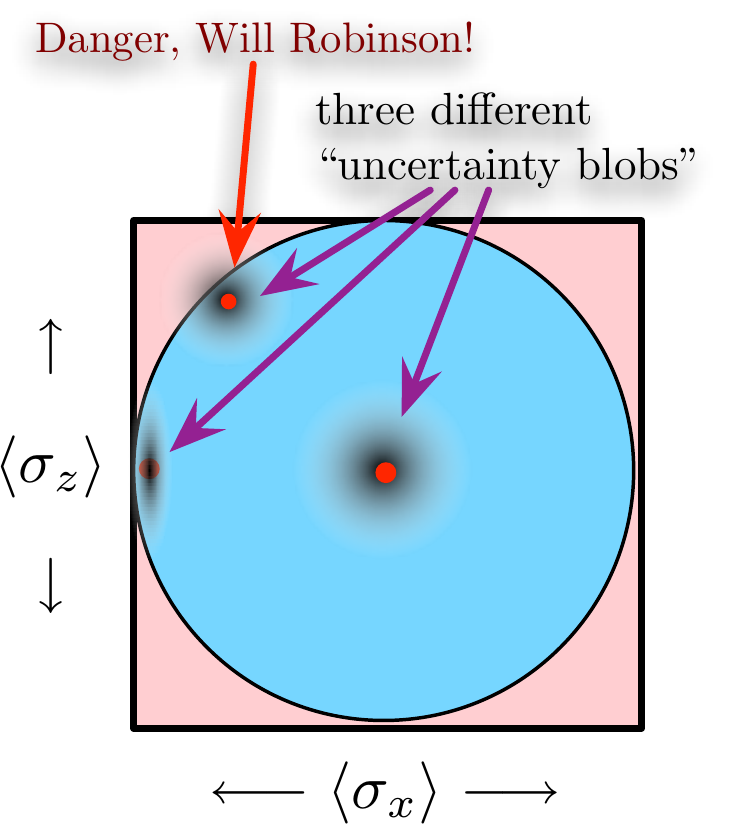}
\caption{The state space of a \emph{qubit} comprises all $2\times2$ positive semidefinite, trace-1 \emph{density matrices}, and is usually represented as the \emph{Bloch ball} (shown here in cross-section).  Quantum state tomography of a single qubit normally proceeds by measuring three 2-outcome observables, $\sigma_x$, $\sigma_z$, and $\sigma_y$ (not shown).  Each can be seen as an independent ``coin'', except that their biases are constrained (by positivity of the density matrix) to the Bloch ball.  When the density matrix is nearly pure (i.e., close to the surface of the Bloch ball), and \emph{not} diagonal in one of the measured bases, the qubit behaves much like a noisy coin.  As shown, linear inversion can easily yield estimates outside the convex set of physical states (analogous to the $p$-simplex).}
\label{fig:QuantumStateEstimation}
\end{figure}

The central idea of quantum information science \cite{NielsenChuang} is to encode information into pure quantum states, then process it using unitary quantum gates, and thus to accomplish certain desirable tasks (e.g., fast factoring of integers, or secure cryptography).

One essential step in achieving these goals is the experimental characterization of quantum hardware using \emph{state tomography} \cite{Tomography} and \emph{process tomography} (which is mathematically isomorphic to state tomography \cite{AncillaAssisted}).

Quantum tomography is, in principle, closely related to probability estimation (see Fig. \ref{fig:QuantumStateEstimation}).  The catch is that quantum states cannot be sampled directly.  Characterizing a single \emph{qubit}'s state requires measuring three independent 1-bit observables.  Typically, these are the Pauli spin operators, $\{\sigma_x,\sigma_y,\sigma_z\}$, which behave like coin flips with respective biases $q_x,q_y,q_z$.  Furthermore, even if the quantum state is pure (i.e., it has zero entropy, and \emph{some} measurement is perfectly predictable), no more than one of these ``coins'' can be deterministic -- and in almost all cases, all three of $q_x,q_y,q_z$ are different from 0 and 1.  So, just as in the case of the noisy coin, we have \emph{pure} states (analogous to $p=0$) that yield somewhat random measurement outcomes ($q\neq0$).  The details are somewhat different from the noisy coin, but sampling mismatch remains the essential ingredient.

The implications of our noisy coin analysis for quantum tomography are more qualitative than in the applications above.  Quantum states are not quite noisy coins, but they are \emph{like} noisy coins.  In particular, the worst-case risk must scale as $O(1/\sqrt{N})$ rather than as $O(1/N)$, for any possible fixed choice of measurements.  Finding the exact minimax risk will require a dedicated study, but the analysis here strongly suggests that HML (first proposed for precisely this problem in \cite{BlumeKohout2010Hedged}) will perform well.

Even more interesting is the implication that \emph{adaptive} tomography can offer an enormous reduction in risk.  This occurs because for quantum states -- unlike noisy coins -- the amount of ``noise'' is under the experimenter's control.  By measuring in the eigenbasis of the true state $\rho$, sampling mismatch can be eliminated!  However, this requires knowing (or guessing) the eigenbasis.  Bagan et al \cite{Bagan2006Separable} observed this, via a somewhat different analysis, and also demonstrated that in the $N\to\infty$ limit it is sufficient to adapt \emph{once}.  An extension of the analysis presented here should be able to determine near-minimax adaptive strategies for finite $N$.

\noindent\textbf{Acknowledgments:}  CF acknowledges financial support from the Government of Canada through NSERC.  RBK acknowledges financial support from the Government of Canada through the Perimeter Institute, and from the LANL LDRD program.

\bibliographystyle{plain}
%\bibliography{csferrie}

\end{document}